\documentclass{amsart}
\usepackage{amsmath, amsthm, amsfonts, amsopn}
        \theoremstyle{plain}
        \newtheorem{Thm}{Theorem}[section]
        \newtheorem{Prop}{Proposition}[section]
        \newtheorem{Lemma}{Lemma}[section]

        \theoremstyle{definition}
        \newtheorem{Def}{Definition}[section]
        \newtheorem{Example}{Example}[section]
        
        \theoremstyle{remark}
        \newtheorem*{Remark}{Remark}
        
        \numberwithin{equation}{section}
        
        \newcommand{\field}[1]{\mathbb{#1}}
        
        \newcommand{\Z}{\field{Z}}
        
        \newcommand{\R}{\field{R}}
        \newcommand{\C}{\field{C}}
        \newcommand{\abs}[1]{\lvert#1\rvert}
        \newcommand{\norm}[1]{\lVert#1\rVert}

        \newcommand{\define}[1]{\emph{#1}}      

        \newcommand{\cover}[1]{\widetilde{#1}}
        
    \newcommand{\Laplace}{\Delta}
    \newcommand{\ltwo}{l^{2}}
    \newcommand{\Ltwo}{L^{2}}
    \DeclareMathOperator{\Tr}{Tr}
    \newcommand{\VNTr}[1]{\Tr_{#1}}
    \DeclareMathOperator{\Det}{Det}
    \newcommand{\FKdet}[1]{\Det_{#1}}
    \newcommand{\mFKdet}[1]{\Det^{\prime}_{#1}}
    \newcommand{\Betti}[2]{b_{#1}(#2)}
    \newcommand{\LtwoBetti}[2]{b_{#1}^{(2)}(#2)}
    \newcommand{\Luck}{L\"{u}ck}

        \newcommand{\Folner}{F\o lner}
        \newcommand{\trivgp}{\{e\}}
        \newcommand{\present}[2]{\left\langle#1\ |\ #2\right\rangle}
        \DeclareMathOperator{\Wh}{Wh}

\begin{document}

\title{Residual Amenability and the Approximation of $\Ltwo$-invariants}

\author{Bryan Clair}
\address{Department Of Mathematics,
         University Of Chicago.
         Chicago, IL 60615.}
\email{bryan@math.uchicago.edu}
\date{October 3, 1997}

\begin{abstract}
We generalize \Luck's Theorem to show that the $\Ltwo$-Betti numbers
of a residually amenable
covering space are the limit of the $\Ltwo$-Betti
numbers of a sequence of amenable covering spaces.
We show that
any residually amenable covering space of a finite simplicial
complex is of determinant class, and that the $\Ltwo$ torsion
is a homotopy invariant for such spaces.
We give examples of residually amenable groups, including the
Baumslag-Solitar groups.
\end{abstract}

\maketitle

\section*{Introduction}

In 1994, Wolfgang \Luck~\cite{luck:resfin} proved the beautiful theorem
that if $X$ is a finite simplicial complex with residually finite
fundamental group, the $\Ltwo$-Betti numbers of the universal covering
of $X$ can be approximated by the ordinary Betti numbers
of a sequence of finite coverings of $X$.
In fact, the question of approximation dates back to
Kazhdan~\cite{kaz:arith} (see also~\cite[Pg. 20]{gr:asympt})
but only an inequality was known.
Dodziuk and Mathai~\cite{dm:amen} have shown a result analagous
to \Luck's Theorem in the situation where the covering transformation
group is amenable.
Specifically, they show that the $\Ltwo$-Betti numbers of an amenable
covering $\cover{X}$ of $X$ can be approximated by the ordinary Betti numbers
of a sequence of \Folner\ subsets of $\cover{X}$.
This paper generalizes \Luck's
Theorem to the case where the cover of $X$ has residually amenable
transformation group,
a large class of groups that includes the residually finite
groups of \Luck's Theorem and the amenable groups of Dodziuk and Mathai.

In this paper, we also consider $\Ltwo$ torsion.
For $\Ltwo$ acyclic covering spaces,
$\Ltwo$ analytic torsion was first studied in~\cite{mat:l2antor}
and~\cite{lott:heat}, and $\Ltwo$ Reidemeister-Franz torsion was first 
studied in~\cite{cm:l2tor}, see also~\cite{lur:kthe}.  
To define these $\Ltwo$ torsions, one needs to establish decay of 
the $\Ltwo$ spectral density function at 0.  In the case of a 
residually finite covering, \Luck~\cite{luck:resfin},
derives an elegant estimate on the spectral distribution functions
for the finite covers, which in the limit gives the necessary decay for the 
combinatorial $\Ltwo$ Laplacian.  \Luck\ also proves the homotopy 
invariance of $\Ltwo$ combinatorial torsion in this case.

Recently, in~\cite{cfm:detlines}, the combinatorial and analytic
torsion invariants were 
defined more generally as volume forms on $\Ltwo$ homology and $\Ltwo$ 
cohomology respectively, the decay condition on the spectrum now replaced by a 
similar condition known as determinant class.
This allowed interpretation of results of~\cite{bfkm:tor} as the 
equality of the combinatorial and analytic $\Ltwo$ torsions.

Dodziuk and Mathai~\cite{dm:amen} show that coverings with amenable 
covering group are of determinant class, and Mathai and
Rothenberg~\cite{mr:homotor} have recently extended \Luck's results to prove 
the homotopy invariance of $\Ltwo$ torsion in
that case.  In this paper we show that coverings with residually 
amenable covering group are of determinant class, and that $\Ltwo$ 
torsion of such spaces is a homotopy invariant.

On a different note, Farber~\cite{far:growth} has recently
generalized \Luck's Theorem in a
new direction, viewing it as a statement about flat bundles
rather than finite coverings.  In particular, he gives precise 
conditions for the convergence of $\Ltwo$-Betti numbers of finite
non-regular covers.  A reasonable direction for future work would be
to try and extend the results of this paper using his techniques.

We now formulate the main results of this paper.
Let $Y$ be a connected simplicial complex.  Suppose that a
finitely generated group $\pi$ acts
freely and simplicially on $Y$ so that $X = Y/\pi$ is a finite
simplicial complex.

Suppose there is a nested sequence of normal subgroups
\newcommand{\grouptower}{
        {\pi = \Gamma_{1} \supset \Gamma_{2} \supset \dotsb}}
$\grouptower$ such that
$\bigcap_{n=1}^{\infty}\Gamma_{n} = \trivgp$.
Form $Y_{n} = Y/\Gamma_{n}$, so that $Y_{1},Y_{2}, \dotsc$ are a
tower of covering spaces of $X$.

Say that $\pi$ is residually finite if there exist $\Gamma_{n}$'s
so that the quotients $\pi/\Gamma_{n}$ are all finite.
Then \Luck's Theorem~\cite{luck:resfin} states that
\[ \LtwoBetti{j}{Y:\pi} =
   \lim_{n \rightarrow \infty} 
              \frac{1}{\abs{\Gamma_{n}}}\Betti{j}{Y_{n}} \]
where $\LtwoBetti{j}{Y:\pi}$ is the $j$th $\Ltwo$-Betti number
of $Y$.

We generalize \Luck's Theorem to the situation
where $\pi$ is residually amenable, meaning there exist $\Gamma_{n}$'s
so that the quotients $\pi/\Gamma_{n}$ are all amenable.
The first main result of this paper is

\begin{Thm}[Approximation Theorem]  \label{thm:mainbetti}
Suppose $Y$ is a simplicial complex, $\pi$ acts freely and
simplicially on $Y$, and $X = Y/\pi$ is a finite simplicial complex.
If $\pi$ is residually amenable, then
\[ \LtwoBetti{j}{Y:\pi} =
   \lim_{n \rightarrow \infty} \LtwoBetti{j}{Y_{n}:\pi/\Gamma_{n}}. \]
\end{Thm}

The next result gives more evidence for the determinant
class conjecture, which states that any regular covering space of a
finite simplicial complex is of determinant class.  For $\pi$ residually
finite this follows from~\cite{luck:resfin}, and it was
shown for $\pi$ amenable in~\cite{dm:amen}.

\begin{Thm}[Determinant Class Theorem]  \label{thm:maindet}
Suppose $Y$ is a simplicial complex, $\pi$ acts freely and
simplicially on $Y$, and $X = Y/\pi$ is a finite simplicial complex.
If $\pi$ is residually amenable, then $Y$ is of determinant class.
\end{Thm}

Now we turn to the problem of homotopy invariance of $\Ltwo$ torsion.  
Let $M$ and $N$ be compact manifolds, $\cover{M}$ and 
$\cover{N}$ regular $\pi$-covering spaces.
As in~\cite{mr:homotor}, a homotopy equivalence $f:M \rightarrow N$
induces a canonical isomorphism
\( \cover{f}^{*}_{*}:\det \overline{H}^{*}_{(2)}(\cover{N})
             \rightarrow \det \overline{H}^{*}_{(2)}(\cover{M}) \)
of determinant lines of $\Ltwo$ cohomology.

Let $\phi_{\cover{M}} \in \det \overline{H}^{*}_{(2)}(\cover{M})$
denote the $\Ltwo$ torsion of $M$.
\begin{Thm}[Homotopy Invariance of $\Ltwo$ Torsion] \label{thm:mainhomo}
Suppose $f:M \rightarrow N$ is a homotopy equivalence of compact odd 
dimensional manifolds, and $\cover{M}$ and 
$\cover{N}$ are regular covering spaces with residually amenable 
covering group.  Then via the above identification of determinant lines of 
$\Ltwo$ cohomology,
\[ \phi_{\cover{M}} = \phi_{\cover{N}} \in
       \det \overline{H}^{*}_{(2)}(\cover{M}). \]
\end{Thm}

This provides more evidence for the conjecture in~\cite{mr:homotor}
(see also~\cite{luck:resfin}) that $\Ltwo$ torsion is 
always a homotopy invariant when the covering spaces in question are 
of determinant class.

This paper is organized as follows.  
The first section covers preliminaries on residually amenable groups,
and exhibits interesting examples.
The second section proves the main technical theorem. It 
essentially states that the $\Ltwo$-spectra of Laplacians on the 
$Y_{n}$ approximate the $\Ltwo$-spectrum of Laplacian on $Y$. 
Finally, the third section proves the three main theorems.

Thanks to Mel Rothenberg for suggesting this direction of research,
and also to Kevin Whyte and Shmuel Weinberger for helpful discussions.

\section{Preliminaries}
\subsection{Residual Properties}
\newcommand{\Class}{\mathcal{C}}

\begin{Def}
Let $\Class$ be a nonempty class of groups (though possibly containing
only one group).
A group $\pi$ is \define{residually $\Class$} if for any element
$g \in \pi$, $g \neq e$, there exists a quotient group $\pi'(g)$ 
belonging to $\Class$ such that $g \mapsto g' \in \pi'(g)$ with
$g' \neq e$.
\end{Def}

If a countable group $\pi$ is residually $\Class$, then there is
a nested sequence of normal subgroups $\grouptower$
such that $\pi/\Gamma_{n}$ belongs to $\Class$ and
$\bigcap_{n=1}^{\infty}\Gamma_{n} = \trivgp$.
For other basic theorems concerning residual properties
of groups, we refer to~\cite{magnus:resfin}.

When $\Class = \{\text{finite groups}\}$ we say that $\pi$ is a 
\define{residually finite} group.

\subsection{Amenability}
Let $\pi$ be a finitely generated discrete group, with word metric $d$. 
We use the following characterization of amenability, due to \Folner.

\begin{Def}
$\pi$ is \define{amenable} if there is a sequence of finite subsets
$\{\Lambda_{k}\}_{k=1}^{\infty}$ such that for any fixed $\delta>0$
\[ \lim_{k \rightarrow \infty}
        \frac{\#\{\partial_{\delta}\Lambda_{k}\}}
                 {\#\{\Lambda_{k}\}}
        = 0 \]
where
\( \partial_{\delta}\Lambda_{k} =
   \{\gamma \in \pi : d(\gamma,\Lambda_{k}) < \delta
                \text{ and }  d(\gamma,\pi - \Lambda_{k}) < \delta\} \)
is a $\delta$-neighborhood of the boundary of $\Lambda_{k}$.
\end{Def}

Examples of amenable groups include finite groups, abelian groups,
nilpotent groups and solvable groups, and groups of subexponential
growth.  Amenability for discrete groups is preserved by the following 
five processes:
\newcounter{x} \begin{list}{\arabic{x}.}
        {\usecounter{x}
                \setlength{\itemsep}{0pt}\setlength{\parsep}{0pt}}
\item Taking subgroups;
\item Forming quotient groups;
\item Forming group extensions by amenable groups;
\item Forming upward directed unions of amenable groups;
\item Forming a direct limit of amenable groups.
\end{list}

Free groups with two or more generators,
and fundamental groups of closed negatively curved 
manifolds are \emph{not} amenable.

\subsection{Residual Amenability}
\begin{Def}
If $\pi$ is residually $\Class$, where $\Class = \{\text{amenable groups}\}$,
we say that $\pi$ is \define{residually amenable}.
\end{Def}

Recall that the derived subroups $\pi^{(i)}$ of a group $\pi$ are defined
by $\pi^{(0)} = \pi$ and $\pi^{(i+1)} = [\pi^{(i)},\pi^{(i)}]$.  Say
that $\pi$ is \define{solvable} if $\pi^{(i)} = \trivgp$ for some $i$.

Contained in the class of residually amenable groups is the important
class of residually solvable groups.  Unlike residually 
finite groups, the class of residually solvable
groups is closed under extensions.  And in 
contrast to amenable groups, the free product of two
residually solvable groups is residually solvable.  The statement for 
extensions is Prop~\ref{prop:ext} below.  Closure under free products
follows from the fact that solvability is a root property as discussed
in~\cite{magnus:resfin}.

\begin{Prop}
The following are equivalent:
\begin{list}{\emph{(\roman{x})}}
        {\usecounter{x}
                \setlength{\itemsep}{0pt}\setlength{\parsep}{0pt}}
\item $\pi$ is residually solvable;
\item $\bigcap_{i=1}^{\infty}\pi^{(i)} = \trivgp$;
\item $\pi$ contains no nontrivial perfect subgroup.  $\Gamma$ is
        \define{perfect} if $\Gamma = [\Gamma,\Gamma]$.
\end{list}
\end{Prop}
\begin{proof}
(ii)$\implies$(i) is clear.  Since $\bigcap_{i=1}^{\infty}\pi^{(i)}$
is a perfect subgroup of $\pi$, (iii) implies (ii).

Now suppose $\pi$ is residually solvable, and suppose that
$\Gamma$ is a perfect subgroup of $\pi$.  Then $\Gamma$ is also
residually solvable.  If $\Gamma$ is nontrivial, there is some
nontrivial map $f : \Gamma \rightarrow S$ with $S$ solvable of rank $k$.
Then 
\[ \trivgp = S^{(k)} \supseteq f(\Gamma)^{(k)} = 
                f(\Gamma^{(k)}) = f(\Gamma) \]
which contradicts nontriviality of $f$.  This shows (i)$\implies$(iii).
\end{proof}

\begin{Prop}\label{prop:ext}
If $\Gamma_{1}$ and $\Gamma_{2}$ are residually solvable, and $\pi$ is
an extension
\[ 1 \rightarrow \Gamma_{1} \stackrel{\iota}{\rightarrow} \pi
      \stackrel{\kappa}{\rightarrow} \Gamma_{2} \rightarrow 1 \]
then $\pi$ is residually solvable.
\end{Prop}
\begin{proof}
Suppose $H$ is a perfect subgroup of $\pi$.  Then
$\kappa(H)$ is a perfect supgroup of $\Gamma_{2}$, hence trivial.
Then $H \subset \Gamma_{1}$ and so $H$ is trivial.  Thus $\pi$ is
residually solvable.
\end{proof}

\begin{Example}
For nonzero integers $p$ and $q$,
define the \define{Baumslag-Solitar} group
$BS(p,q)$ by
\[ BS(p,q) = \present{a,b}{a^{-1}b^{p}a = b^{q}}. \]

A group $\pi$ is \define{Hopfian} if $\pi/\Gamma \cong \pi$
implies $\Gamma = \trivgp$.
The family of groups $BS(p,q)$ were first defined in~\cite{bs:bs},
where it was shown that $BS(p,q)$ is Hopfian if and only if $p$ and $q$
are \define{meshed}, which means $p|q$, $q|p$, or $p$ and $q$ have exactly
the same set of prime divisors.
As any finitely generated residually finite
group is Hopfian, the groups $BS(p,q)$ are not residually finite when
$p$ and $q$ are not meshed.

Kropholler shows in~\cite{krop:bs} that the second derived subgroup
$\pi^{(2)}$ is free when $\pi = BS(p,q)$, for any $p$, $q$.
When $p$ and $q$ are not both $\pm 1$, $\pi^{(2)}$ is free on two
or more generators and therefore $\pi$ is not amenable.
However, $\pi$ is residually solvable since it is an extension
\[
        1 \rightarrow \pi^{(2)} \rightarrow \pi \rightarrow
        \pi/\pi^{(2)} \rightarrow 1
\]
with $\pi^{(2)}$ residually solvable and $\pi/\pi^{(2)}$ solvable.

More generally, let $\pi$ be any non-cyclic group which is the fundamental
group of a graph of infinite cyclic groups. Then from~\cite{krop:bs},
$\pi^{(2)}$ is free and the above argument shows $\pi$ is residually
solvable.
\end{Example}

\begin{Example}
We show that the amalgamated free product of two abelian groups is 
residually solvable.  Slightly more generally, suppose we have two
residually solvable groups $A$ and $B$, and two homomorphisms from
a group $H$ into the centers of $A$ and $B$ given by
$\alpha : H \rightarrow Z(A)$ and $\beta : H \rightarrow Z(B)$.
Then the free product with amalgamation $A*_{H}B$ is residually 
solvable.

To see this, let \( N = \{(\alpha(h),\beta(h)^{-1})|h \in H\}
\subset A \times B. \)  Since $H$ is abelian, $N$ is a subgroup of
$A \times B$. $N$ is normal because $H$ includes into the centers
of both $A$ and $B$.  Note that if $H$ is only known to be normal
in $A$ and $B$, $N$ is unlikely to be normal in $A \times B$.

Let $K$ be the kernel of the natural map
$A*_{H}B \rightarrow (A \times B)/N$.  In $A*_{H}B$, $K$ has trivial 
intersection with all conjugates of $A$ and $B$, hence $K$ is free by 
a well known theorem of group actions on trees.  Then $A*_{H}B$ is
an extension of the residually solvable group $(A \times B)/N$ by
the residually solvable group $K$, and so $A*_{H}B$ is residually 
solvable.
\end{Example}

\begin{Example}
It is shown in~\cite{rv:reshnn} that any HNN-extension of a finitely
generated abelian group is residually solvable.
\end{Example}

\begin{Example}
R.J. Thompson, G. Higman, K. Brown, and E.A. Scott have demonstrated
various classes of finitely presented infinite simple groups.
(See for example,~\cite{scott:infsimple}).  The example of Scott
contains a free subgroup on two generators and is therefore not
amenable and not residually amenable.
\end{Example}

\section{Main Technical Theorems}
%
%

As in the introduction, suppose $Y$ is a simplicial complex,
$\pi$ acts freely and
simplicially on $Y$, and $X = Y/\pi$ is a finite simplicial complex.
Suppose we have a nested sequence of normal
subgroups $\grouptower$ such that
$\bigcap_{n=1}^{\infty}\Gamma_{n} = \trivgp$.
Define $Y_{n} = Y/\Gamma_{n}$.

        \newcommand{\fund}{\mathcal{F}}                                 
        
Suppose $X$ has $a_{j}$ cells in dimension $j$, and choose
a lift to $Y$ of each $j$-cell of $X$.  These choices give a basis
over $\ltwo(\pi)$
of the space $C^{j}_{(2)}(Y)$ of $j$ dimensional $\ltwo$-cochains on $Y$.
The lifts also descend to give bases
of $C^{j}_{(2)}(Y_{n})$ over $\ltwo(\pi/\Gamma_{n})$.

Denote by $\Laplace$ and ${\Laplace}_{n}$ the Laplacian on
$C^{j}_{(2)}(Y)$ and $C^{j}_{(2)}(Y_{n})$ respectively.  All
arguments to follow will apply to a specific value of $j$, but
this dependence will not be indicated.

Let $\{P(\lambda) : \lambda \in [0,\infty)\}$ and
$\{P_{n}(\lambda) : \lambda \in [0,\infty)\}$ denote the right
continuous family of spectral projections of $\Laplace$ and 
${\Laplace}_{n}$.  
Since $\Laplace$ is $\pi$-equivariant,  so are $P(\lambda) = 
\chi_{[0,\lambda]}(\Laplace)$ for $\lambda \in [0,\infty)$.  
Similarly, $P_{n}(\lambda)$ are $\pi/\Gamma_{n}$-equivariant.  Let
$F, F_{n} : [0,\infty) \rightarrow [0,\infty)$ denote the spectral density
functions
               \[ F(\lambda) = \VNTr{\pi}P(\lambda) \]
               \[ F_{n}(\lambda) = \VNTr{\pi/\Gamma_{n}}P(\lambda). \]

We now set
\[ \begin{array}{rclrcl}
      \overline{F}(\lambda) & = &
                \limsup_{n \rightarrow \infty} F_{n}(\lambda); &
      \underline{F}(\lambda) & = &
                \liminf_{n \rightarrow \infty} F_{n}(\lambda) \\
      \overline{F}^{+}(\lambda) & = &
                \lim_{\delta \rightarrow 0^{+}}\overline{F}(\lambda+\delta);&
      \underline{F}^{+}(\lambda) & = &
                \lim_{\delta \rightarrow 0^{+}}\underline{F}(\lambda+\delta).\\
      
\end{array} \]

With the above notation, we state the main technical results
of this paper.

\begin{Thm}  \label{thm:tech1}
For all $\lambda \in [0,\infty)$,
        \[ F(\lambda) = \overline{F}^{+}(\lambda)
                                                = \underline{F}^{+}(\lambda)  \]
\end{Thm}

\begin{Thm}  \label{thm:tech2}
Suppose there is some right continuous function
$s : [0,\varepsilon) \rightarrow [0,\infty)$,
with $s(0) = 0$ and
so that for all $n$ and for all $\lambda \in [0,\varepsilon)$ we have
\[ F_{n}(\lambda) - F_{n}(0) \leq s(\lambda) \]
then
\begin{enumerate}
        \item $\overline{F}(\lambda)$ and $\underline{F}(\lambda)$ are
                right continuous at zero and one has the equalities
                \[      \overline{F}(0) = \overline{F}^{+}(0) = F(0) =
                        \underline{F}(0) = \underline{F}^{+}(0). \]
        \item For all $\lambda \in [0,\varepsilon)$,
                                \[ F(\lambda) - F(0) \leq s(\lambda). \]
\end{enumerate}
\end{Thm}

These theorems and their proofs are quite similar to
\Luck~\cite[Theorem 2.3]{luck:resfin}, but here they are
stated so as to require no conditions
on the quotient groups $\pi/\Gamma_{n}$.  The residual finiteness
condition in \Luck's Theorem or the residual amenability condition
of this paper are required to provide $s$, the uniform decay at zero 
of the spectral density functions for the covers $Y_{n}$.

To show the two technical theorems, we first prove a
number of preliminary lemmas.

\begin{Lemma}\label{lem:norm}
There exists a number $K > 1$ such that the operator norms of
$\Laplace$ and $\Laplace_{n}$ are smaller than $K^{2}$
for all $n$.
\end{Lemma}
\begin{proof}
Choosing lifts of cells of $X$, we have identified the space of
$\ltwo$-cochains on $Y$ with $\bigoplus_{i=1}^{a}\ltwo(\pi)$.
The combinatorial Laplacian $\Laplace$ is then described by an
$a \times a$ matrix $B$ with entries in $\Z[\pi]$, acting by
right multiplication.  The Laplacian $\Laplace_{n}$ is described by 
the same matrix $B$, now acting by right multiplication on
$\bigoplus_{i=1}^{a}\ltwo(\pi/\Gamma_{n})$

For $u = \sum_{g \in \pi}\lambda_{g}\cdot g \in \C[\pi]$ define
$\abs{u}_{1} = \sum_{g \in \pi}\abs{\lambda_{g}}$.  Choose
$K > 1$ so that
\[ K \geq a \cdot \sum_{j=1}^{b}
                  \max{\left\{\abs{B_{ij}}_{1}; i=1,\dotsc,a\right\}}. \]
The proof then proceeds exactly as in~\cite[Lemma 2.5]{luck:resfin}.
\end{proof}

\begin{Lemma}\label{lem:res}
Let $p(\mu)$ be a polynomial.  There is a number $n_{0}$, depending 
only on the system of groups \( \pi = \Gamma_{1} \supset \Gamma_{2} 
\supset \ldots \) such that for all $n \geq n_{0}$
\[  \VNTr{\pi}p(\Laplace) = \VNTr{\pi/\Gamma_{n}}p(\Laplace_{n}) \]
\end{Lemma}
\begin{proof}
We identify $\Laplace$ with an
$a \times a$ matrix $B$ with entries in $\Z[\pi]$, as in the
previous Lemma.

Fix elements $g_{0}, g_{1},\ldots,g_{r} \in \pi$ and
$\lambda_{0}, \lambda_{1},\ldots,\lambda_{r} \in \R$
such that $g_{0} = e$, $g_{i} \neq e$, and $\lambda_{i} \neq 0$ for
$1 \leq i \leq r$ so that
\[ \sum_{j=1}^{a}(p(B))_{j,j} = \sum_{i=0}^{r}\lambda_{i}g_{i}. \]
Then \[ \VNTr{\pi}p(\Laplace) = \lambda_{0} .\]

The Laplacian $\Laplace_{n}$ on $Y_{n}$ is also described by the
matrix $B$, now acting on $\bigoplus_{i=1}^{a}\ltwo(\pi/\Gamma_{n})$
by right multiplication.  Then
\[ \VNTr{\pi/\Gamma_{n}}p(\Laplace_{n}) =
   \sum_{i=1}^{r}\lambda_{i}\VNTr{\pi/\Gamma_{n}}R(g_{i}) \]
where
$R(g_{i}) : \ltwo(\pi/\Gamma_{n}) \rightarrow \ltwo(\pi/\Gamma_{n})$
is right multiplication with $g_{i}$.

Since the intersection of the $\Gamma_{i}$'s is trivial, there is a
number $n_{0}$ such that for $n \geq n_{0}$ none of the elements
$g_{i}$ for $1 \leq i \leq r$ lies in $\Gamma_{n}$.  Since $\Gamma_{n}$
is normal, we conclude for $n \geq n_{0}$ and $i \neq 0$
\[ \VNTr{\pi/\Gamma_{n}}R(g_{i}) = 0 .\]

Then for $n \geq n_{0}$
\[ \VNTr{\pi}p(\Laplace) = \lambda(0) =
   \VNTr{\pi/\Gamma_{n}}p(\Laplace_{n}). \]

\end{proof}

\begin{Lemma}\label{lem:pchi}
Let $\{p_{k}(\mu)\}_{k=1}^{\infty}$
be a sequence of polynomials, uniformly bounded on
$[0,\norm{\Laplace}]$, such that for the
characteristic function $\chi_{[0,\lambda]}(\mu)$ of the interval
$[0,\lambda]$,
\[ \lim_{k \rightarrow \infty}p_{k}(\mu) = \chi_{[0,\lambda]}(\mu) \]
holds for each $\mu \in [0,\norm{\Laplace}]$.  Then
\[ \lim_{k \rightarrow \infty}\VNTr{\pi}p_{k}(\Laplace) = 
   F(\lambda). \]
\end{Lemma}
\begin{proof}
        This lemma and its proof are identical to \cite[Lemma 2.7]{luck:resfin}
\end{proof}

        \newcommand{\onek}{\textstyle{\frac{1}{k}}}             

We now prove Theorem~\ref{thm:tech1}.  Fix $\lambda \geq 0$.  Define
for $k \geq 1$ a continuous function $f_{k} : \R \rightarrow \R$ by
\[ f_{k}(\mu) = 
        \begin{cases}
                1 + \onek                                       & \mu \leq \lambda \\
                1 + \onek - k(\mu - \lambda)& \lambda \leq \mu \leq \lambda + \onek \\
                \onek                                           & \lambda + \onek \leq \mu
        \end{cases}
\]
Clearly $\chi_{[0,\lambda]}(\mu) < f_{k+1}(\mu) < f_{k}(\mu)$ and
$f_{k}(\mu)$ converges to 
$\chi_{[0,\lambda]}(\mu)$ for all $\mu \in [0,\infty)$.  For each $k$
choose a polynomial $p_{k}$ such that
$\chi_{[0,\lambda]}(\mu) < p_{k}(\mu) < f_{k}(\mu)$ holds for all
$\mu \in [0,K^{2}]$, where $K$ is as in Lemma~\ref{lem:norm}.
The polynomials $p_{k}$ satisfy the conditions of Lemma~\ref{lem:pchi}.

Because $\chi_{[0,\lambda]}(\mu) \leq p_{k}(\mu)$ for all $\mu \in
[0,\norm{\Laplace_{n}}]$, we have
\begin{equation} \label{e:lbd} \begin{split}
        F_{n}(\lambda)& =
                        \VNTr{\pi/\Gamma_{n}}\bigl(\chi_{[0,\lambda]}(\Laplace_{n})\bigr) \\
        & \leq  \VNTr{\pi/\Gamma_{n}}\bigl(p_{k}(\Laplace_{n})\bigr).
\end{split} \end{equation}

On the other hand, we have $p_{k}(\mu) \leq 1+\onek$
for $\mu \in [0, \lambda + \onek]$ and $p_{k}(\mu) \leq \onek$
for $\mu \in [\lambda + \onek, K^{2}]$.  So
\begin{equation}\label{e:ubd}\begin{split}
        \VNTr{\pi/\Gamma_{n}}\bigl(p_{k}(\Laplace_{n}) \bigr) &\leq
                \VNTr{\pi/\Gamma_{n}}
                \bigl((1+\onek)\chi_{\left[0,\lambda + \onek\right]}(\Laplace_{n}) 
                  \bigr) \\
        & \quad + \VNTr{\pi/\Gamma_{n}}
                \bigl((\onek)\chi_{\left[\lambda + \onek,K^{2}\right]}(\Laplace_{n})
                  \bigr) \\
        &= (1 + \onek)F_{n}(\lambda + \onek) \\
        & \quad + \onek(F_{n}(K^{2}) - F_{n}(\lambda + \onek)) \\
        &= F_{n}(\lambda + \onek) + \onek F_{n}(K^{2})
\end{split} \end{equation}

Now notice 
\(
   F_{n}(K^{2}) = 
   \VNTr{\pi/\Gamma_{n}}\bigl(\chi_{[0,K^{2}]}(\Laplace_{n})\bigr).
\)
But $\chi_{[0,K^{2}]}(\Laplace_{n})$ is the identity on the space
$C^{j}_{(2)}(Y_{n})$, which is identified with
$\bigoplus_{i=1}^{a_{j}}\ltwo(\pi/\Gamma_{n})$.
Thus
\begin{equation}\label{e:unif}
        F_{n}(K^{2}) = a_{j}.
\end{equation}

By Lemma~\ref{lem:res}, there is a number $n_{0}(k)$ for each
polynomial $p_{k}$ such that for $n \geq n_{0}(k)$
\[  \VNTr{\pi}p_{k}(\Laplace) = 
        \VNTr{\pi/\Gamma_{n}}\bigl(p_{k}(\Laplace_{n})\bigr). \]
Then for $n \geq n_{0}(k)$, the equations \eqref{e:lbd}, \eqref{e:ubd}, 
and \eqref{e:unif} give
\[ F_{n}(\lambda) \leq
   \VNTr{\pi}p_{k}(\Laplace) \leq
   F_{n}(\lambda + \onek) + \onek a_{j} \]
Taking $lim_{n \rightarrow \infty}$,
\[ \overline{F}(\lambda) \leq
   \VNTr{\pi}p_{k}(\Laplace) \leq
   \underline{F}(\lambda + \onek) + \onek a_{j} \]
Taking $lim_{k \rightarrow \infty}$, and using Lemma~\ref{lem:pchi},
\[ \overline{F}(\lambda) \leq
   F(\lambda) \leq
   \underline{F}^{+}(\lambda). \]
We have for all $\epsilon > 0$
\[ F(\lambda) \leq \underline{F}^{+}(\lambda)
    \leq \underline{F}(\lambda + \epsilon)
    \leq \overline{F}(\lambda + \epsilon)
    \leq F(\lambda + \epsilon) \]
and since
$\lim_{\epsilon \rightarrow 0+}F(\lambda + \epsilon) = F(\lambda)$
we get
\[ F(\lambda) = \overline{F}^{+}(\lambda)
   = \underline{F}^{+}(\lambda). \]

This finishes the proof of Theorem~\ref{thm:tech1}. \qed

Next we show Theorem~\ref{thm:tech2}.1 and Theorem~\ref{thm:tech2}.2.
We suppose there is some right continuous function 
$s : [0,\varepsilon) \rightarrow [0,\infty)$,
with $s(0) = 0$
and
so that for all $n$ and for all $\lambda \in [0,\varepsilon)$ we have
\[ F_{n}(\lambda) - F_{n}(0) \leq s(\lambda). \]

Taking the limit inferior and limit superior for $n \rightarrow 
\infty$ gives:
\[ \underline{F}(\lambda) \leq \underline{F}(0) + s(\lambda)
\text{ and }
   \overline{F}(\lambda) \leq \overline{F}(0) + s(\lambda). \]
Taking the limit for $\lambda \rightarrow 0$ gives
\[ \underline{F}^{+}(0) \leq \underline{F}(0) \text{ and }
   \overline{F}^{+}(0) \leq \overline{F}(0). \]
And finally, since $\underline{F}$ and $\overline{F}$ are increasing,
\[ \underline{F}^{+}(0) = \underline{F}(0) \text{ and }
   \overline{F}^{+}(0) = \overline{F}(0). \]

We already know $\overline{F}^{+}(0) = F(0) = \underline{F}^{+}(0)$ from
Theorem~\ref{thm:tech1}, and this proves Theorem~\ref{thm:tech2}.1.
Since $s$ is right continuous, we conclude:
\[ \overline{F}^{+}(\lambda) \leq F(0) + s(\lambda) \]
and Theorem~\ref{thm:tech2}.2 follows from Theorem~\ref{thm:tech1}.
This finishes the proof of Theorem~\ref{thm:tech2}.
\qed

The following lemma is needed in the proof of 
Theorem~\ref{thm:maindet} in the last section.

        \newcommand{\specint}[1]{{\int_{0^{+}}^{K^{2}}
                        \frac{F_{#1}(\lambda)-F_{#1}(0)}{\lambda}d\lambda}}

\begin{Lemma}\label{lem:int}
\[ \specint{}  \leq  \liminf_{n \rightarrow \infty} \specint{n} \]
\end{Lemma}
\begin{proof}
The proof of this lemma is identical to that of Lemma 3.3.1 in
\cite{luck:resfin}.  It follows from Theorem~\ref{thm:tech1} and
the monotone convergence theorem.
\end{proof}

%
%
\section{Proofs Of The Main Theorems}
In this section, we prove the approximation theorem and determinant
class theorem for residually amenable groups.

\subsection{The Approximation Theorem}
\begin{proof}[Proof of Theorem~\ref{thm:mainbetti} (Approximation 
Theorem)]~   

Observe that the $j$th $\Ltwo$ Betti numbers of $Y$ and $Y_{n}$ are
given by
     \[ \begin{array}{cc}
           b^{j}_{(2)}(Y:\pi) = F(0); &
           b^{j}_{(2)}(Y_{n}:\pi/\Gamma_{n}) = F_{n}(0).
         \end{array} \]
Therefore Theorem~\ref{thm:mainbetti} will follow directly
from Theorem~\ref{thm:tech2}.1, if we can establish a uniform decay of
$F_{n}$ near zero.

Since $\pi/\Gamma_{n}$ is amenable,~\cite[Theorem 2.1.3]{dm:amen}
applies, and we have a constant $K > 1$ so that
\begin{equation}\label{eq:decay}
        F_{n}(\lambda) - F_{n}(0) \leq 
        a_{j}\frac{\log K^{2}}{-\log \lambda} = s(\lambda)
\end{equation}
for all $0 < \lambda < 1$.  The constant $K$ can be
any number larger than $\max(\norm{\Laplace_{n}},1)$, and
therefore can be chosen independently of $n$, by Lemma~\ref{lem:norm}.
\end{proof}

\begin{Remark}
In the case of \Luck's Theorem, the groups $\pi/\Gamma_{n}$ are
finite, so the Laplacian on $Y_{n}$ is a finite, self-adjoint matrix.
\Luck~ then proves, for any self-adjoint matrix, an elegant estimate on
the number of eigenvalues which are
less than a fixed $\lambda$.  The estimate is weakend if the product
of the nonzero eigenvalues is small, but this product must be at
least one as the Laplacian has integer entries.

Dodziuk and Mathai make use of the same fundamental estimate when
proving the theorem for amenable groups.
\end{Remark}

\subsection{Results For Manifolds}
Suppose now $M$ is a compact Riemannian manifold, $\cover{M}$ is a 
regular
covering space for $M$ with residually amenable transformation group,
and let $\cover{\Laplace}$ denote the Laplacian on $\Ltwo\ j$-forms 
on $\cover{M}$.
Following arguments in~\cite{dm:amen},
one can investigate the analytic spectral
density function $\cover{F}(\lambda)$ of $\cover{\Laplace}$.
To relate this analytic Laplacian to the combinatorial situation in 
the previous sections, choose $X$ to be a triangulation of $M$
and lifting to a triangulation $Y$ of $\cover{M}$. 

The results analagous to those of~\cite{dm:amen} are:

\begin{Thm}(Gap Criterion)
The spectrum of $\cover{\Laplace}$ has a gap at zero if and only if
there is a $\lambda > 0$ such that
\[ \lim_{n \rightarrow \infty} F_{n}(\lambda) - F_{n}(0) = 0. \]
\end{Thm}

\begin{Thm}(Spectral Density Estimate)
There are constants $C > 0$ and $\varepsilon > 0$ independent of 
$\lambda$, such that for all $\lambda \in (0,\varepsilon)$
\[ \cover{F}(\lambda) - \cover{F}(0) \leq \frac{C}{-\log{\lambda}}.\]
\end{Thm}

\subsection{The Determinant Class Theorem}
A covering space $Y$ of a finite simplicial complex $X$ is said to be 
of \define{determinant class} if, for each $j$,
\[ -\infty < \int_{0^{+}}^{1}\log\lambda dF(\lambda), \]
where $F(\lambda)$ denotes the von Neumann spectral density function 
of the combinatorial Laplacian $\Laplace_{j}$ on $\Ltwo j$-cochains.

For a regular cover of a compact Riemannian manifold $M$, there is
a corresponding notion of \define{analytic determinant class}.
We can also ask if $\cover{M}$ is of 
determinant class in the combinatorial sense, by choosing a 
triangulation of $M$.  This question turns out to be 
independent of triangulation, and equivalent to analytic determinant 
class.  For a more detailed discussion, see~\cite{dm:amen}.

We will prove that every residually amenable covering of a finite 
simplicial complex is of determinant class.  The appendix 
of~\cite{bfk:gluing} 
contains a proof that every residually finite covering of a compact 
manifold is of determinant class.  Their proof is based on \Luck's 
approximation of von Neumann spectral density functions.  Since an 
analagous approximation holds in the setting of this paper, we can 
apply the argument of~\cite{bfk:gluing} to prove Theorem~\ref{thm:maindet}.
The fact that our coverings are infinite neccessitates some
modifications.

        \newcommand{\modspec}[2]
                {{(\log #2)\bigl(F_{#1}(#2) - F_{#1}(0)\bigr)}}

\begin{proof}[Proof of Theorem~\ref{thm:maindet} (Determinant 
Class Theorem)]~

As with the rest of this paper, this proof will proceed for a fixed 
$j$ which will be supressed in the notation.

Denote by $\mFKdet{\pi/\Gamma_{n}}\Laplace_{n}$ the modified 
Fuglede-Kadison~\cite{fk:det} determinant of $\Laplace_{n}$, that is, the 
Fuglede-Kadison determinant of $\Laplace_{n}$ restricted to the 
orthogonal complement of its kernel.  It is given by the following 
Lebesgue-Stieltjes integral,
\begin{equation*}
   \log \mFKdet{\pi/\Gamma_{n}}\Laplace_{n} =
   \int_{0^{+}}^{K^{2}}\log\lambda dF_{n}(\lambda)
\end{equation*}
with $K$ as in Lemma~\ref{lem:norm}.  That is, $\norm{\Laplace_{n}} 
\leq K^{2}$.

Integration by parts yields
\begin{equation}\label{eq:parts}
   \begin{split}
    \log \mFKdet{\pi/\Gamma_{n}}\Laplace_{n}&
    = \modspec{n}{K^{2}} \\
    & + \lim_{\epsilon \rightarrow 0^{+}}
      \left\{
          -\modspec{n}{\epsilon} 
              -{\int_{\epsilon}^{K^{2}}
                                        \frac{F_{n}(\lambda)-F_{n}(0)}{\lambda}d\lambda}
          \right\}.
   \end{split}
\end{equation}

Now $Y_{n}$ is an amenable cover of $X$, so we are in the
situation of Dodziuk-Mathai~\cite{dm:amen}.  We require two results
contained in their proof
of the Determinant Class Theorem for amenable 
coverings~\cite[Thm 0.2]{dm:amen}.
The first is that
\[ \log \mFKdet{\pi/\Gamma_{n}}\Laplace_{n} \geq 0, \]
which is stronger than simply determinant class.
The second is the existence of the integral
\begin{equation}\label{eq:specint} \specint{n}. \end{equation}

It is worth remarking that their proof and
the corresponding proof in~\cite{bfk:gluing}
require similar statements for the Laplacians on finite 
approximations.  In the finite case, the Laplacian is an integer 
matrix, so the determinant is a positive integer and therefore at 
least 1.  The existence of the integral (\ref{eq:specint}) is clear
in the finite case, as the spectrum is discrete.

From (\ref{eq:parts}) and the existence of (\ref{eq:specint}),
$\lim_{\epsilon \rightarrow 0^{+}}\modspec{n}{\epsilon}$
exists.  In fact, for the integral (\ref{eq:specint}) to exist, we 
must have
\[
        \lim_{\epsilon \rightarrow 0^{+}}\modspec{n}{\epsilon} = 0.
\]

Then (\ref{eq:parts}) becomes
\begin{equation*}
        \log \mFKdet{\pi/\Gamma_{n}}\Laplace_{n}
    = \modspec{n}{K^{2}} - \specint{n}
\end{equation*}
which gives
\begin{equation*}
        \specint{n} \leq \modspec{n}{K^{2}}.
\end{equation*}

Now from Lemma~\ref{lem:int},
\begin{equation}\label{eq:bound}
\begin{split}
        \specint{} &\leq \liminf_{n \rightarrow \infty} \specint{n} \\
                           &\leq \liminf_{n \rightarrow \infty} \modspec{n}{K^{2}}.
\end{split}
\end{equation}
Since we have uniform decay~(\ref{eq:decay}) of the $F_{n}$ near 0,
Theorem~\ref{thm:tech2}.1 applies and
$\lim_{n \rightarrow \infty}F_{n}(0) = F(0)$.  Since
$K^{2} \geq \norm{\Laplace_{n}}$ for all $n$
and $K^{2} \geq \norm{\Laplace}$,
\[ F_{n}(K^{2}) = F(K^{2}) = a_{j} \]
for all $n$, so (\ref{eq:bound}) becomes
\begin{equation}\label{eq:shoot}
        \specint{} \leq \modspec{}{K^{2}}.
\end{equation}

This shows in particular that the left hand integral exists, and
arguing as with $\Laplace_{n}$ we have
\[
        \log \mFKdet{\pi}\Laplace = \modspec{}{K^{2}} - \specint{}
\]
which is non-negative by~(\ref{eq:shoot}).
Since this is true for all $j$, $Y$ is of determinant class.

\end{proof}                     

\subsection{Homotopy Invariance of $\Ltwo$ Torsion}
The Fuglede-Kadison determinant $\FKdet{\pi}$ induces a homomorphism
from the Whitehead group of $\pi$,
\[  \Phi_{\pi}:\Wh(\pi) \rightarrow \R^{+} \]
which was defined in~\cite{lur:kthe} and~\cite{luck:resfin}.
Given a homotopy eqivalence $f:M \rightarrow N$ of compact manifolds, 
we choose cell decompositions for $M$ and $N$, choose $f$ to be a 
cellular homotopy equivalence, and let $M_{f}$ be the cellular mapping 
cone.  Then the cochain complex $C^{*}(M_{f})$ is an acyclic complex 
over the group ring $\Z[\pi]$, and defines the Whitehead torsion 
$T(f) \in \Wh(\pi)$.  Then $\Phi_{\pi}(T(f)) \in \R^{+}$.

Now suppose $\pi$ is residually amenable, so $M$ and $N$ are of 
determinant class.
Let $\phi_{\cover{M}}$, $\phi_{\cover{N}}$
denote the $\Ltwo$ torsion of $M$ and $N$.  The homotopy equivalence
$f$ canonically identifies the determinant lines of $\Ltwo$ 
cohomology of $\cover{M}$ and $\cover{N}$, so
$\phi_{\cover{M}}\otimes\phi_{\cover{N}}^{-1} \in \R^{+}$.

Then from~\cite[Prop. 2.1]{mr:homotor}, one has
\[ \phi_{\cover{M}}\otimes\phi_{\cover{N}}^{-1} =
     \Phi_{\pi}(T(f)) \in \R^{+}. \]
Therefore, Theorem~\ref{thm:mainhomo} will follow from the following
\begin{Prop}
Suppose that $\pi$ is a finitely presented residually amenable group.
Then the homomorphism
\[ \Phi_{\pi}:\Wh(\pi) \rightarrow \R^{+} \]
is trivial.
\end{Prop}

\begin{proof}

\newcommand{\LapKL}[1]{\Laplace^{\cover{K}_{#1},\cover{L}_{#1}}}
We can represent an arbitrary element of $\Wh(\pi)$ as the Whitehead 
torsion of a homotopy equivalence $f:L\rightarrow K$ of finite 
CW-complexes, which without loss of generality is an inclusion.
Let $\cover{L}$ and $\cover{K}$ denote the corresponding regular 
$\pi$ covering complexes.  The relative cochain complex 
$C^{*}(\cover{K},\cover{L})$ is acyclic, and so is its $\Ltwo$ 
completion $C^{*}_{(2)}(\cover{K},\cover{L})$.  In particular, the 
combinatorial Laplacian $\Laplace_{j}^{\cover{K},\cover{L}}$ is 
invertible and we see that
\[ \Phi_{\pi}(T(f)) =
        \prod_{j=0}^{n}
                {\FKdet{\pi}(\LapKL{}_{j})}^{\frac{(-1)^{j}j}{2}}
        > 0.
\]
We claim that $\FKdet{\pi}(\LapKL{}_{j}) = 1$ for
each $j$.  The $j$ will be supressed in the notation.

Form the covering spaces $\cover{L}_{n}$ and $\cover{K}_{n}$ with 
amenable covering group $\pi/\Gamma_{n}$, and let
$\LapKL{n}$ denote the Laplacian on
$C^{*}_{(2)}(\cover{K}_{n},\cover{L}_{n})$.
As $\LapKL{}$ is invertible we can choose a 
common bound $K$ on $\norm{\LapKL{}}$ and $\norm{(\LapKL{})^{-1}}$,
as in Lemma~\ref{lem:norm}.
$K$ will also bound $\norm{\LapKL{n}}$ and $\norm{(\LapKL{n})^{-1}}$.

It then follows from the proof of~\cite[Prop. 2.6]{mr:homotor} that
for all $n$,
\[ \FKdet{\pi/\Gamma_{n}}\LapKL{n} = 1 \]
and $\LapKL{n}$ has a spectral gap at zero of size at least $K^{-2}$.
Then, analagous to~\cite[Thm 3.4.3]{luck:resfin},
\[ \FKdet{\pi}\LapKL{} =
        \lim_{n \rightarrow \infty}\FKdet{\pi/\Gamma_{n}}\LapKL{n} = 1. \]
\end{proof}
        
This finishes the proof of Theorem~\ref{thm:mainhomo}.
As remarked earlier, in both~\cite{luck:resfin} and ~\cite{mr:homotor},
the approximating spaces are finite.  The approximating Laplacians
have determinant 1 and spectral gap because they are finite
invertible integer matrices, and the previous proposition 
essentially boils down to this fact.

%
%
\bibliographystyle{abbrv}
\bibliography{bryanrefs}

\end{document}